\documentclass[12pt]{amsart}

\usepackage{amsfonts}
\usepackage{amssymb}
\usepackage{graphicx}
\usepackage{amsmath,mathtools}
\usepackage{latexsym}
\usepackage{amscd}
\usepackage{xypic}
\usepackage{mathrsfs}%scrfont}
\usepackage{enumitem} %http://ctan.org/pkg/enumitem
\usepackage{braket}
\usepackage[margin=3cm]{geometry}
\usepackage{amsthm}
\usepackage{accents}
\usepackage{xcolor}

\usepackage{hyperref}

\setlist{itemsep=3pt}

\allowdisplaybreaks

\numberwithin{equation}{section}

\newtheorem{prop}{Proposition}
\newtheorem{theo}[prop]{Theorem}%changed
\newtheorem*{theo*}{Theorem}%changed
\newtheorem{lemm}[prop]{Lemma} %changed
\newtheorem{coro}[prop]{Corollary} %changed

\theoremstyle{definition}

\newtheorem{rema}[prop]{Remark}
\newtheorem{defi}[prop]{Definition}
\newtheorem*{defi*}{Definition}

\numberwithin{prop}{section}

\newcommand{\p}{\partial}
\newcommand{\BB}{\mathbb{B}}

\newcommand{\RR}{\mathbb{R}} %changed
\renewcommand{\SS}{\mathbb{S}}
\newcommand{\TT}{\mathbb{T}}

\newcommand{\cA}{\mathcal A}

\newcommand{\cC}{\mathcal C}

\newcommand{\cE}{\mathcal E}

\newcommand{\cG}{\mathcal G}
\renewcommand{\cH}{\mathcal H}

\renewcommand{\cR}{\mathcal R}

\newcommand{\cV}{\mathcal V}

\DeclareMathOperator{\Ric}{Ric}

\newcommand{\sff}{\mathrm{I\!I}}

\usepackage{graphicx}
\allowdisplaybreaks

\newenvironment{nouppercase}{%
  \renewcommand{\uppercasenonmath}[1]{}}{}

\title{Comparison Theorems and the Intermediate Ricci Curvature Assumption }

\author{Yujie Wu}
\begin{document}
\begin{nouppercase}
	\maketitle
\end{nouppercase}

\begin{abstract}
We explore the notion of $m$-intermediate Ricci curvature assumption introduced by Brendle-Hirsch-Johne further. If a manifold has non-negative $m$-intermediate Ricci curvature and stable weighted slicing of order $m-1$, $\Sigma_0\supset...\supset \Sigma_{m-1}$, then the last slice $\Sigma_{m-1}$ has almost non-negative Ricci curvature in the spectral sense. We prove comparison theorems on the diameter and in-radius bound for stable weighted (respectively free boundary) slicing in such manifolds (respectively with mean convex boundary).
\end{abstract}

\section{Introduction}

In \cite{brendle2024generalization}, Brendle-Hirsch-Johne introduced the notion of ``$m$-intermediate curvature'' on a Riemannian manifold $M^n$ ($1\leq m\leq n$), where $1$-intermediate curvature is equal to Ricci curvature, and $n$-intermediate Ricci curvature is equal to a multiple of scalar curvature.

For a Riemannian manifold $X$, given $p\in X$ and $e_p,e_q\in T_p X$ we denote the sectional curvature as $K_X(e_p,e_q):=\text{Rm}_X(e_p,e_q,e_q,e_p)$ where $\text{Rm}(\cdot,\cdot,\cdot,\cdot)$ is the Riemann curvature tensor with the convention 
$$\text{Rm}(e_p,e_q,e_q,e_p)=\langle\nabla_{e_p}\nabla_{e_q}e_q-\nabla_{e_q}\nabla_{e_p}e_q-\nabla_{[e_p,e_q]}e_q,e_p\rangle.$$

\begin{defi}[Brendle-Hirsh-Johne, \cite{brendle2024generalization}] Suppose $(X^n,g)$ is a Riemannian manifold (possibly with boundary). Given a point $p\in X$ and  orthonormal vectors $\{e_1,...,e_m\}$ of $T_pX$, we define the $m$-intermediate Ricci curvature,
\begin{equation*}C_m(e_1,...,e_m):=\sum_{p=1}^m\sum_{q=p+1}^n K_X(e_p,e_q)
\end{equation*} 
	We say that $(X,g)$ has non-negative $m$-intermediate curvature at $p\in X$ if $C_m(e_1,...,e_m)\geq 0$ for any choice of orthonormal vectors $\{e_1,...,e_m\} \subset T_p X$, and we say the manifold $X$ has non-negative $m$-intermediate curvature if it has non-negative $m$-intermediate curvature at any point $p\in X$.
\end{defi}

\begin{rema}
	Having $C_2 \geq 0$ is the same as having non-negative  Bi-Ricci curvature as introduced by Shen and Ye in \cite{shen1996stable}.
\end{rema}

It's a classical theorem of Bonnet and Meyers that $M^m\times \SS^1$ does not admit a metric of positive Ricci curvature for any closed manifold $M$; and Geroch's conjecture states that the torus $\TT^n$ has no complete metric of positive scalar curvature, this was resolved by Schoen and Yau using minimal hypersurfaces for $n\leq 7$ \cite{yau1979structure} and by Gromov and Lawson for all dimensions using spinors \cite{gromov1983positive}. Furthermore, Chodosh and Li generalized this conjecture in proving that aspherical manifolds of dimension $4$ or $5$ also do not admit a metric of positive scalar curvature \cite{chodosh2024generalized}.
 Brendle-Hirsch-Johne (\cite{brendle2024generalization}) proved that this can be generalized to an obstruction of $m$-intermediate curvature for manifolds of the type $M^{n-m}\times \TT^m$.

%We define the following constant as appeared in \cite{brendle2024generalization}, for $n\geq 2, 1\leq m\leq n-1,k=n-m$,
%
%\begin{equation*}
%	C(m,n):=\frac{m^2-2-n(m-2)}{2(n-m)(m-1)}=\frac{m^2-2-(m+k)(m-2)}{2k(m-1)}.
%\end{equation*}

\begin{theo}[Brendle-Hirsh-Johne, \cite{brendle2024generalization}]\label{BHJ}
	Let $(N^n,g)$ be a closed Riemannian manifold, $2\leq m\leq n-1, m^2-2-n(m-2)\geq 0$ and $N$ has $C_m>0$, then $N$ has no stable weighted slicing of order $m$ (such that the last slice is compact). If $n\leq 7$ 
	(in this case we always have $m^2-2-n(m-2)\geq 0$), suppose $(N^n,g)$ is a closed orientable Riemannian manifold with a degree-nonzero map $F: N^n\rightarrow M^{n-m}\times \TT^m$ for some closed orientable manifold $M$, then $N$ has a stable weighted slicing of order $m$; in particular, any such manifold $N$ does not admit a metric of positive $m$-intermediate curvature.
\end{theo}
We will recall the definition of ``stable weighted slicing of order $m$'' in section \ref{prelim}.

Theorem \ref{BHJ} has been extended by Chen (\cite{chen2024generalization}) to manifolds of the type $M^{n-m}\times \TT^m\times X$ for a closed manifold $M$ and an arbitrary manifold $X$ when $3\leq n\leq 5$, or when $n-m\leq 2, 6\leq n\leq 7$.
On the other hand, manifolds of the type $M^{n-m}\times \TT^m$ can have non-negative $m$-intermediate curvature. In this case, Chu-Kwong-Lee (\cite{chu2025rigidity}) proved that when $n\leq 5$ and $\cC_m \geq 0$, then $M^{n-m} \times \TT^m$ must be isometrically covered by  $X^{n-m} \times \RR^m$ for some compact $X$ with $\Ric_X \geq 0$; Xu (\cite{xu2025dimension}) later extended this to the case when $n=6$ and provided counterexamples to Theorem \ref{BHJ} when $m^2-2-n(m-2)<0$.

In this paper we look at manifolds of the type $N^n=M^{k+1}\times \TT^{m-1}$ with $n=k+m$, where we can find stable minimal slicing of order $m-1$, and if $\cC_m>0$, then the last slice $\Sigma_{m-1}$ have almost positive Ricci curvature in the spectral sense, that is $\Sigma^{k+1}:=\Sigma_{m-1}$ is a manifold with a smooth positive function $f$, such that for some $\gamma>0,\alpha>0$,
\begin{align*}
	-\gamma\Delta_{\Sigma}f+\Ric_{\Sigma}f\geq -\alpha |\nabla_{\Sigma}\log f|^2.
\end{align*}

For suitable choices of $\gamma, \alpha,k,m$, we can prove comparison theorems for such slicings, generalizing the results of Shen-Ye on Bi-Ricci curvature.

\begin{theo}[Shen-Ye, \cite{shen1996stable}]
	Assume that $\cC_2\geq k>0$ and $m\leq 5$, then any closed stable minimal hypersurface $S^{m-1}$ in $M^m$ has $\text{diam}(S)\leq \frac{\pi}{\sqrt{c(m)k}}$, where $c(3)=\frac{3}{4}, c(4)=\frac{1}{2}, c(5)=\frac{1}{4}$.
\end{theo}

Generalization of the result of Shen-Ye has been applied to the slicing of $\mu$-bubbles, in obtaining obstructions to positive intermediate curvature assumptions in aspherical manifolds under suitable dimensional restrictions by Mazurowski-Wang-Yao (\cite{mazurowski2025topology}). The role of the dimensional restriction to such obstructions has been studied by Xu, in showing that there is a complete metric on $\SS^4 \times \RR^2$ with uniformly positive $\cC_2$ and infinite Urysohn width (\cite{xu2025dimension}). Below we proved a sharp control in diameter for such slicing.
\begin{theo}[Bonnet-Meyers Theorem for the $(m-1)$-th Slice]\label{BM-Cm}
	Take a complete Riemannian manifold $N^n$ with positive intermediate curvature $C_m\geq k_0>0$, let $n<m+2+\frac{2}{m-1}$ (this implies $m^2-2-n(m-2)\geq 0$). If we have stable minimal slicing of order $m-1$ for $\max\{1,n-6\}\leq m\leq n-1$,  $N=\Sigma_0 \supset \Sigma_1 \supset \Sigma_2... \supset \Sigma_{m-1}$, then $\text{diam}(\Sigma_{m-1})\leq \frac{\pi}{\sqrt{k_0C_0}}$, where
	$$C_0:=-\frac{n(m-1)-m(m+1)}{2(m^2-2-n(m-2))}>0, \quad \text{when } n<m+2+\frac{2}{m-1}.$$
	
	The constant for $m=2$ is the same constant as obtained in \cite{shen1996stable}. The bound is sharp and obtained by the embedding of $\SS^3 \hookrightarrow \TT^{m-1} \times \SS^3$, with $C_m \geq 2$ and $C_0=\frac{1}{2}$.
\end{theo}
%\begin{coro}
%	If $N=M^{k+1}\times \TT^{m-1}$ has $\cC_m>k_0>0$, then $\text{diam}(M)\leq \frac{\pi}{\sqrt{k_0C_0}}$ when $k<2+\frac{2}{m-1}$ and $1\leq m\leq \min\{n-1,7\}$.
%\end{coro}
\begin{rema}
	We note that the compactness of $\Sigma_{m-1}$ is not assumed a priori.
\end{rema}

%We note the only solutions of $C(m,n)=0$ when $m,n\in \NN$ are when $n=7$ and $m=3$ or $m=4$. \textcolor{red}{Are there any counterexample for $n=6,m=2$? or $m=3, n=6$ (the latter case still has $C(m,n)>0$ but does not have $C_0>0$. Any other case $C_0>0$ and $C(m,n)>0$ are equivalent.)}

Assuming $C_m\geq 0$ and $m$-convexity, we have the following theorem for stable free boundary weighted slicing of order $m$. 

Let $(X,\p X,g)$ be a Riemannian manifold with boundary, and we denote the second fundamental form $\sff_{\p X}(Y,Z)=-\langle \nabla_Y Z,\nu_{\p X} \rangle$ for $\nu_{\p X}$ the outward pointing unit normal along $\p X$ and $Y,Z$ vector fields along $\p X$.
\begin{defi}[$m$-convexity assumption]
	 Take $(X^n,\p X,g)$ a Riemannian manifold with boundary, $1\leq m\leq n-1$, we say $X$ is (strictly) $m$-convex at a point $p\in X$ if for any choice of orthonormal vectors $\{e_1,...,e_m\}\subset T_p {\p X}$, $\sum_{i=1}^m\sff_{\p X}(e_i,e_i)>0$. We denote $\sff^{\p X}_m\geq k_0>0$ if $\sum_{i=1}^m\sff_{\p X}(e_i,e_i)>k_0$ for any choice of orthonormal vectors at any point on $\p X$ and say $X$ is uniformly $m$-convex.
\end{defi}
\begin{theo}[Stable Weighted Free Boundary Slicing] \label{NoFBSlicing}
	Let $(X^n,\p X,g)$ be a Riemannian manifold with boundary,  assume $C_m\geq 0$ and $\sff^{\p X}_m >0$, if $m^2-2-n(m-2)\geq 0$ then $X$ has no stable free boundary weighted slicing of order $m$ such that the last slice $\Sigma_m$ is compact.
\end{theo}

%Example: $\BB^k \times \TT^m$ has boundary equal to $\SS^{k-1} \times \TT^m$ does not have $\sff_m>0$.

\begin{coro}
\label{RelHomology}
	Let $Y$ be a compact orientable smooth manifold, then any  smooth orientable Riemannian manifold $(N^n, \p N)$ that maps to $(Y,\p Y)\times \TT^m$ with non-zero degree ($\p N$ is mapped to $\p Y \times \TT^m$), has a stable free boundary weighted slicing of order $m$ when $n\leq 7$. If $m^2-2-n(m-2)\geq 0$, the product manifold $(Y,\p Y)\times \TT^m$ has no metric of non-negative intermediate curvature and uniformly $m$-convex boundary. 
\end{coro}

For stable weighted free boundary slicing we have the following comparison result for the in-radius (see \cite{gromov2019mean}).
\begin{theo}[In-Radius Bound]\label{InRadiusTh}
	If a manifold $X^n$ has intermediate curvature $C_m \geq 0$ and uniformly positive mean curvature $H_{\p X}\geq H_0 >0$. Let let $n<m+2+\frac{2}{m-1}$ (this implies $m^2-2-n(m-2)\geq 0$), if we have stable minimal free boundary slicing of order $m-1$ for $\max\{1,n-6\}\leq m\leq n-1$,  $X=\Sigma_0\supset \Sigma_1 ...\supset \Sigma_{m-1}$, and the last slice has compact $\p \Sigma_{m-1}$, then the in-radius of $\Sigma_{m-1}$ is uniformly bounded, $\text{Rad}_{in}(\Sigma_{m-1}):=\sup_{x\in \Sigma_{m-1}}d_{\Sigma_{m-1}}(x,\p \Sigma_{m-1})\leq \frac{1}{C_0H_0}$, for the same $C_0$ defined in Theorem \ref{BM-Cm}.
\end{theo}

Note that $\BB^3 \times \TT^{m-1}$ has $C_0=\frac{1}{2}$, $\cC_m \geq 0, H_0=2$ and every point in $\BB^3\times \{z\}$ for $z\in \TT^{m-1}$ is at most distance $1$ from the boundary.

\textbf{Organization of the Article.} In section \ref{prelim}, we recall the notion of stable weighted slicing and relevant computations involving $\mu$-bubbles. In section \ref{BM-slice}, we prove Theorem \ref{BM-Cm} which generalizes the result of Shen-Ye (\cite{shen1996stable}) on stable minimal hypersurfaces to stable weighted slicing of order $m-1$ on a manifold of positive $m$-intermediate curvature. In section \ref{BM-FBslice-InRad}, we prove Theorem \ref{NoFBSlicing} and Corollary \ref{RelHomology}, extending the result of Brendle-Hirsh-Johne (\cite{brendle2024generalization}) to stable weighted free boundary slicing; we also prove Theorem \ref{InRadiusTh} which extends the In-Radius bound for manifolds with non-negative Ricci curvature and uniformly mean convex boundary, to manifolds with non-negative intermediate Ricci curvature and uniformly mean convex boundary. Lastly, in section \ref{BMSpec}, we prove analogous comparison theorems as in section \ref{BM-slice} and section \ref{BM-FBslice-InRad},  for manifolds with non-negative intermediate curvature and uniformly positive mean curvature in the spectral sense.

\textbf{Acknowledgement.} 
	{Funded by the European Union (ERC Starting Grant 101116001 – COMSCAL). Views and opinions expressed are however those of the author(s) only and do not necessarily reflect those of the European Union or the European Research Council. Neither the European Union nor the granting authority can be held responsible for them.}

\section{Preliminary}\label{prelim}

We first recall the notion of stable weighted slicing introduced in \cite{brendle2024generalization}, an idea that appeared in Schoen and Yau's resolution of Geroch conjecture in \cite{yau1979structure}. 

\begin{defi}[Stable weighted slicing of order $m$, \cite{brendle2024generalization}]Let $1\leq m\leq n-1$ and $(N^n,g)$ a Riemannian manifold. A stable weighted slicing of order $m$ consists of a collection of two-sided and smooth immersions $\Sigma_k\hookrightarrow \Sigma_{k-1}, 1\leq k\leq m$ and a collection of positive functions $\rho_k \in C^{\infty}(\Sigma_k)$ such that,
\begin{itemize}
	\item $\Sigma_0=N, \rho_0=1$;
	\item For each $1\leq k\leq m$, $\Sigma_k$ is a two-sided hypersurface in $\Sigma_{k-1}$ and a stable critical point to the weighted area functional among hypersurfaces in $\Sigma_{k-1}$,
	\begin{equation*}\cH_k(\Sigma):=\int_{\Sigma} \rho_{k-1} d\cH^{n-k}
	\end{equation*}
	\item For $1\leq k\leq m$, $\rho_k=v_k\cdot\rho_{k-1}\rvert_{\Sigma_k}$ for $v_k$ a positive first eigenfunction of the stability operator associated to $\cH^k$.
\end{itemize} 
\end{defi}

We have the following corresponding defintion for free boundary slicing.
\begin{defi}[Stable free boundary weighted slicing of order $m$] \label{DefFBSlice}
	Let $1\leq m\leq n-1$ and $(X^n,\p X, g)$ a Riemannian manifold with boundary. A stable free boundary weighted slicing of order $m$ consists of a collection of two-sided and smooth immersions $(\Sigma_k,\p \Sigma_k)\hookrightarrow (\Sigma_{k-1},\p \Sigma_{k-1}), 1\leq k\leq m$ and a collection of positive functions $\rho_k \in C^{\infty}(\Sigma_k)$ such that,
\begin{itemize}
	\item $\Sigma_0=X, \rho_0=1$
	\item For each $1\leq k\leq m$, $\Sigma_k$ is a two-sided free boundary hypersurface in $\Sigma_{k-1}$ and a stable critical point to the weighted area functional among variations whose restriction to $\p \Sigma_{k-1}$ lies in $T\p\Sigma_{k-1}$ (variations that maps $\p \Sigma_{k-1}$ to $\p \Sigma_{k-1}$),
	\begin{equation*}\cH_k(\Sigma):=\int_{\Sigma} \rho_{k-1} d\cH^{n-k}
	\end{equation*}
	\item For $1\leq k\leq m$, $\rho_k=v_k\cdot\rho_{k-1}\rvert_{\Sigma_k}$ for $v_k$ a positive first eigenfunction of the stability operator associated to $\cH^k$.
\end{itemize}
\end{defi}

Our approach is the method of $\mu$-bubbles introduced by Gromov, a generalization of soap bubbles (minimal hypersurfaces).

\begin{defi}Given a Riemannian manifold $(N^n,g)$ and functions $\rho\in C^{\infty}(N)$ and $h \in C^{\infty}(U)$ where $U$ is an open subset of $N$; fix an open set $\Omega_0 \subset U$ with smooth non-empty boundary $\p^* \Omega_0$ in $U$. A (warped) $\mu$-bubble $\Omega$ is a Caccioppoli set in $U$, $\Omega \bigtriangleup\Omega_0 \Subset U$ that is  a stable critical point to the following functional,
\begin{equation*}\cA_h(\Omega)=\int_{\Sigma}\rho d\cH^{n-1}+\int_{N}\rho h(\chi_{\Omega}-\chi_{\Omega_0})d\cH^n,
\end{equation*}
where we denote $\Sigma=\p^*\Omega$  the reduced boundary of a Caccioppoli set (see \cite{maggi2012sets} Chapter 15 for the definition).
\end{defi}

We recall the first and second variation formula for $\mu$-bubbles (which also gives the second variation for the stable weighted slicing). A reference for the proof can be found in \cite{chodosh2024generalized} or \cite{brendle2024generalization}.

\begin{lemm}\label{SecVar}
	If the reduced boundary $\Sigma=\p^* \Omega$ is smooth  with $\Omega$ a $\mu$-bubble to the functional $\cA_h$ as defined above, then we have over $\Sigma$,
	\begin{align*}
		h=&H_{\Sigma}+\nabla_{\nu_{\Sigma}}\log \rho=H_{\Sigma}+\frac{\nabla_{\nu_{\Sigma}}\rho}{\rho},\\
		\frac{d^2}{ds^2}\bigg\rvert_{s=0}=&\int_{\Sigma}\rho f^2(-f^{-1}\Delta_{\Sigma} f-(\Ric_{N}(\nu_{\Sigma},\nu_{\Sigma})+|\sff_{\Sigma}|^2))\\
		&+\int_{\Sigma}\rho f^2(\nabla_{\Sigma}^2\log\rho(\nu_{\Sigma},\nu_{\Sigma})-\langle  \nabla_{\Sigma}\log \rho,\nabla_{\Sigma}\log f\rangle-\nabla_{\nu_{\Sigma}}h),
	\end{align*}
	for any variation with normal speed equal to $f\in C^{\infty}_c(\Sigma)$.
\end{lemm}

Similarly the definition for free boundary $\mu$-bubbles is the following.

\begin{defi}Given a Riemannian manifold with boundary $(N^n, \p N,g)$ and functions $\rho\in C^{\infty}(N)$ and $h \in C^{\infty}(U)$ where $U$ is an open subset of $N$ with $U \cap \p N \neq \emptyset$; fix an open set $\Omega_0 \subset U$ with smooth non-empty boundary $\p^* \Omega_0$ in $U$. A (warped) free boundary $\mu$-bubble $\Omega$ is a Caccioppoli set in $U$, $\Omega \bigtriangleup\Omega_0 \Subset U$ that is  a stable critical point to the following functional,
\begin{equation*}\cA_h(\Omega)=\int_{\Sigma}\rho d\cH^{n-1}+\int_{N}\rho h(\chi_{\Omega}-\chi_{\Omega_0})d\cH^n,
\end{equation*}
\end{defi}

The first and second variation formula for free boundary $\mu$-bubbles are the following. A reference for the proof can be found in \cite{chodosh2024generalized} or \cite{wu2023free}.

\begin{lemm}\label{SecVar}
	If the reduced boundary $\Sigma=\p^* \Omega$ is smooth  with $\Omega$ a free boundary $\mu$-bubble to the functional $\cA_h$ as defined above, then we have that $\Sigma$ meets with $\p N$ orthogonally, and over $\Sigma$, for any variation with normal speed equal to $f\in C^{\infty}_c(\Sigma)$
	\begin{align*}
		h=&H_{\Sigma}+\nabla_{\nu_{\Sigma}}\log \rho=H_{\Sigma}+\frac{\nabla_{\nu_{\Sigma}}\rho}{\rho},\\
		\frac{d^2}{ds^2}\bigg\rvert_{s=0}&=\int_{\p \Sigma} \rho f(\nabla_{\nu_{\p \Sigma}}f-\sff_{\p N}(\nu_{\Sigma},\nu_{\Sigma})f) +
		\int_{\Sigma}\rho f^2(-f^{-1}\Delta_{\Sigma} f-(\Ric_{N}(\nu_{\Sigma},\nu_{\Sigma})+|\sff_{\Sigma}|^2))\\
		&+\int_{\Sigma}\rho f^2(\nabla_{\Sigma}^2\log\rho(\nu_{\Sigma},\nu_{\Sigma})-\langle  \nabla_{\Sigma}\log \rho,\nabla_{\Sigma}\log f\rangle-\nabla_{\nu_{\Sigma}}h),
	\end{align*}
	where $\sff_{\p N}(X,Y)$ denotes $-\langle \nabla_{X}Y,\nu_{\p N}\rangle$ for the outward pointing unit normal $\nu_{\p N}$
\end{lemm}

From now on we denote $\Ric_N(\nu):=\Ric_N(\nu,\nu)$.

\section{Bonnet-Meyers For Slicing of Order $m-1$}\label{BM-slice}

We prove Theorem \ref{BM-Cm} in this section. First we recall a computation in \cite{brendle2024generalization}.

\begin{lemm}[\cite{brendle2024generalization} Lemma 3.1, Lemma 3.2]\label{BHJ1}
	Given stable weighted slicing of order $m-1$, $N=\Sigma_0\supset ... \supset \Sigma_{m-1}$ and the associated weight $\rho_k(1\leq k\leq m-1)$, we have 
	\begin{align*}
		\Delta_{\Sigma_{m-1}}\log\rho_{m-1}&=\sum_{k=1}^{m-1}(-\lambda_k)+\sum_{k=1}^{m-1} (H^2_{\Sigma_k}-|\sff_{\Sigma_k}|^2-\Ric_{\Sigma_{k-1}}(\nu_k))-\sum_{k=1}^{m-1}\langle \nabla_{\Sigma_k}\log \rho_k \nabla_{\Sigma_k} \log v_k \rangle
	\end{align*}
\end{lemm}

\begin{proof}[Proof of Theorem \ref{BM-Cm}]
We may assume $m\geq 2$ as when $m=1$, $C_0=\frac{1}{n-1}$ and Theorem \ref{BM-Cm} is the classical Bonnet-Meyers theorem.

Take any point $p\in N$ and any $0<r<r(p):=\sup_{q\in N}d_{\Sigma_{m-1}}(p,q)$, we want to prove $r \leq C_0$.

	We solve the following minimizing $\mu$-bubble problem for $\cA_{m}$ on the last slice $\Sigma_{m-1}$ for Caccioppoli sets $\Omega$ with $\Sigma=\p^*\Omega$ in $\Sigma_{m-1}$,
	$$\cA_{m}(\Sigma)=\int_{\Sigma}\rho_{m-1}-\int_{\Sigma_{m-1}}\rho_{m-1}h(\chi_{\Omega}-\chi_{\Omega_0}).$$
	
	Here take $h$ to be a smooth mollifier of $\tilde{h}(x):=\beta\tan(\frac{\pi}{r}d_{\Sigma_{m-1}}(x,p)-\frac{\pi}{2}-\epsilon)=:\beta\tan(\phi(x))$ defined over the set $$\tilde{U}:=\{x\in \Sigma_{m-1}, \frac{\epsilon r}{\pi}< d_{\Sigma_{m-1}}(x,p)< r+\frac{\epsilon r}{\pi}<r(p)\},$$
	for some small $\epsilon$ and $\Omega_0=\{x\in \Sigma_{m-1}, h(x)>\epsilon\}$ with $\epsilon$ a regular value of $h$ and some $\beta$ to be chosen later.	
	
	By Proposition 12 in \cite{chodosh2024generalized}, a minimizing $\mu$-bubble $\Sigma_m \rightarrow \Sigma_{m-1}$ must exists and is smooth when the dimension of $\Sigma_{m-1}$ is no more than $7$ (the dimension of $\Sigma_{m-1}$ is equal to $n-(m-1)\leq 7$ by assumption of $m\geq n-6$). So by Lemma \ref{SecVar} we have the following for any compactly supported function $f$ on $\Sigma_m$,
	\begin{align*}
		0\leq &\int_{\Sigma}\rho_{m-1} f^2(-f^{-1}\Delta f-(\Ric_{N}(\nu_m)+|\sff_{\Sigma_m}|^2))\\
		&+\int_{\Sigma}\rho_{m-1} f^2(\nabla^2\log\rho_{m-1}(\nu_m,\nu_m)-\langle  \nabla_{\Sigma_m}\log \rho_{m-1},\nabla_{\Sigma_m}\log f\rangle-\nabla_{\nu_{m}}h),
	\end{align*}
	denoting $\nu_k$ for a chosen unit normal of $\Sigma_k \rightarrow \Sigma_{k-1}$.
	
	Since $\Sigma_m$ is compact by the $\mu$-bubble construction, we can plug-in $f^{-1}=\rho_{m-1}\rvert_{\Sigma_m}$ into the above inequality to get,
	\begin{align*}
		0\leq \int_{\Sigma_{m}}\rho^{-1}_{m-1}\left( \Delta_{\Sigma_{m}}\log \rho_{m-1} -(\Ric_{\Sigma_{m-1}}(\nu_m)+|\sff_{\Sigma_m}|^2)+\nabla^2\log\rho_{m-1}(\nu_m,\nu_m)-\nabla_{\nu_m}h\right)
	\end{align*}
	 
	Using the first variation of $\Sigma_m$ and Lemma 3.1 in \cite{brendle2024generalization}, we have,
	\begin{align*}
		\Delta_{\Sigma_m}\log\rho_{m-1}&=\Delta_{\Sigma_{m-1}} \log \rho_{m-1}+H_{\Sigma_m}(H_{\Sigma_m}-h)-\nabla^2 \log \rho_{m-1}(\nu_m,\nu_m)
	\end{align*}

	Together with Lemma \ref{BHJ1} we get,	
	\begin{align}
		0\leq &\int_{\Sigma_m}\rho^{-1}_{m-1}(-\sum_{k=1}^{m-1}\lambda_k-\sum_{k=1}^m(\Ric_{\Sigma_{k-1}}(\nu_k)+|\sff_{\Sigma_k}|^2-H^2_{\Sigma_k})-\sum_{k=1}^{m-1}\langle \nabla_{\Sigma_k}\log \rho_k \nabla_{\Sigma_k} \log v_k\rangle)\label{3.1}\\
		&-\int_{\Sigma_m} \rho^{-1}_{m-1} (\nabla_{\nu_m}h+H_{\Sigma_m }h)\nonumber\\
		&=-\int_{\Sigma_m} \rho^{-1}_{m-1}(\Lambda+\cR+\cE+\cG +\nabla_{\nu_m}h+H_{\Sigma_m}h),\nonumber
	\end{align}
	the terms $\Lambda, \cR, \cE, \cG$ are as defined in \cite{brendle2024generalization},
	\begin{align*}
		\Lambda&=\sum_{k=1}^{m-1}\lambda_k, \cG=\sum_{k=1}^{m-1}\langle  \nabla_{\Sigma_k}\log \rho_k,\nabla_{\Sigma_k} \log v_k\rangle\\ 
		\cR&=\sum_{k=1}^m \Ric_{\Sigma_{k-1}}(\nu_{k}), \cE=\sum_{k=1}^m|\sff_{\Sigma_k}|^2-H^2_{\Sigma_k}
	\end{align*}
	Applying Lemma 3.10 in \cite{brendle2024generalization} we have,
	\begin{align*}
		\cR+\cE+\cG\geq \cC_m(\nu_1,...,\nu_m)+\sum_{k=1}^m \cV_k.
	\end{align*}
	Lemma 3.11, Lemma 3.12 in \cite{brendle2024generalization} imply that $\cV_k\geq 0$ for $1\leq k\leq m-1$ assuming $C(m,n)\geq 0$,
\begin{equation*}
	C(m,n):=\frac{m^2-2-n(m-2)}{2(n-m)(m-1)}=\frac{m^2-2-(m+k)(m-2)}{2k(m-1)}.
\end{equation*}
	
	On the last slice $\Sigma_m$, different from Lemma 3.7 in \cite{brendle2024generalization}, instead we have 
	\begin{align*}
		\cV_m&=|\sff_{\Sigma_m}|^2- H^2_{\Sigma_m}+\left(\frac{1}{2}+\frac{1}{2(m-1)}\right)(H_{\Sigma_m}-h)^2\\
		&\geq \left(\frac{1}{n-m}-\frac{1}{2}+\frac{1}{2(m-1)}\right)H_{\Sigma_m}^2+\frac{m}{2(m-1)}(-2hH_{\Sigma_m}+h^2)\\
		&=C(m,n)H_{\Sigma_m}^2+\frac{m}{2(m-1)}(-2hH_{\Sigma_m}+h^2)
	\end{align*}
	
Together with $\cC_m \geq k_0>0$, now inequality (\ref{3.1}) becomes 
\begin{align*}
	0\leq -\int_{\Sigma_m}\rho^{-1}_{m-1}(k_0+\nabla_{\nu_m}h-\frac{1}{m-1}hH_{\Sigma_m}+C(m,n)H^2_{\Sigma_m}+\frac{m}{2(m-1)}h^2)
\end{align*}
	We now have,
	\begin{align*}
		&-C(m,n)H^2_{\Sigma_m}+\frac{1}{m-1}hH_{\Sigma_m}\\
		&=-C(m,n)(h^2-2h\nabla_{\nu_m}\log \rho_{m-1}+(\nabla_{\nu_m}\log \rho_{m-1})^2)+\frac{1}{m-1}(h^2-h\nabla_{\nu_m}\log \rho_{m-1})\\
		&=-(C(m,n)-\frac{1}{m-1})h^2+(2C(m,n)-\frac{1}{m-1})h\nabla_{\nu_m}\log \rho_{m-1}-C(m,n)(\nabla_{\nu_m}\log \rho_{m-1})^2\\
		&\stackrel{(\star)}{\leq} -(C(m,n)-\frac{1}{m-1})h^2+\frac{1}{4}(2C(m,n)-\frac{1}{m-1})^2\frac{1}{C(m,n)}h^2\\
		&= \frac{1}{4C(m,n)(m-1)^2}h^2,
	\end{align*}
	where in ($\star$) we used Young's inequality.
	
	Using $$-\frac{m}{2(m-1)}+\frac{1}{4C(m,n)(m-1)^2}=\frac{n(m-1)-m(m+1)}{2(m^2-2-n(m-2))}=:-C_0$$
	Solving $C_0>0$ and $C(m,n)>0$ we get $n<m+2+\frac{2}{m-1}$ for $m\geq 2$.
	
	So we have,
	\begin{align*}
		0\leq & \int_{\Sigma_{m}}-\rho^{-1}_{m-1}(C_0 h^2+\nabla_{\nu_m}h+k_0)\\
		\stackrel{(\star 1)}{\leq} &\int_{\Sigma_m} -\rho^{-1}_{m-1} (C_0\beta^2 \tan^2\phi(x)-\beta \frac{\pi}{r}\cos^2\phi(x)+k_0+\varepsilon) \\
		\stackrel{(\star 2)}{=}& \int_{\Sigma_m} -\rho^{-1}_{m-1}\cdot \left(-(\frac{\pi}{r})^2 C_0^{-1} +k_0+ \varepsilon\right) ,
	\end{align*}
	where in $(\star 1)$ we can choose $\varepsilon>0$ to be arbitrarily small by choosing the mollification $h$ of $\tilde{h}(x)=\beta \tan\phi(x)$; in $(\star 2)$ we let $C_0\beta^2 =\beta \frac{\pi}{r}$ by choosing $\beta$. Now letting $\varepsilon \rightarrow 0$ we obtained $r\leq \frac{\pi}{\sqrt{k_0C_0}}$.

	When $m=2$, $C(m,n)>0, n<m+2+\frac{2}{m-1}=6$, and $C_0=\frac{6-n}{4}$ this agrees with the constant given in Shen and Ye \cite{shen1996stable}.
\end{proof}

\section{stable weighted free boundary slicing and In-Radius Bound}\label{BM-FBslice-InRad}
We now consider existence of stable weighted free boundary minimal slicing $X=\Sigma_0\supset \Sigma_1 \supset ... \supset \Sigma_m$ as defined in Definition \ref{DefFBSlice}. Then by the second variation formula on each $\Sigma_k$ for $1\leq k\leq m$ we have a smooth $v_k>0$,
\begin{align*}
	0\leq \lambda_k=&-\Delta_{\Sigma_{k}} \log v_k-(\Ric_{\Sigma_k-1}(\nu_k)+|\sff_{\Sigma_k}|^2)+\nabla^2_{\Sigma_{k-1}} \log \rho_{k-1} (\nu_k,\nu_k)\\
	&-\langle \nabla_{\Sigma_k}\log \rho_{k-1}, \nabla_{\Sigma_k}\log v_k\rangle-|\nabla_{\Sigma_k} \log v_k|^2,\\
	\nabla_{\nu_{\p \Sigma_k}} \log v_k=& \sff_{\p \Sigma_{k-1}}(\nu_k,\nu_k), \quad \text{along } \p \Sigma_k.
\end{align*}

As in Lemma \ref{BHJ1} we have the following result along the boundary.

\begin{lemm}\label{FB-BC}
	We have,
	\begin{equation*}
		\nabla_{\nu_{\p \Sigma_{k}}}\log {\rho_{k}}=\nabla_{\nu_{\p \Sigma_{k-1}}}\log {\rho_{k-1}}+\sff_{\p \Sigma_{k-1}}(\nu_k,\nu_k)
	\end{equation*}
\end{lemm}

\begin{proof}
	This follows from $\rho_k=\rho_{k-1}v_k$ and $\nu_{\p \Sigma_m}=...=\nu_{\p \Sigma_1}=\nu_{\p\Sigma_0}=\nu_{\p X}$ by free boundary,
	\begin{align*}
		\nabla_{\nu_{\p \Sigma_k}}\log \rho_k=&\nabla_{\nu_{\p \Sigma_k}} \log v_k+\nabla_{\nu_{\p \Sigma_{k-1}}}\log {\rho_{k-1}}\\
		=& \nabla_{\nu_{\p \Sigma_{k-1}}}\log {\rho_{k-1}}+\sff_{\p \Sigma_{k-1}}(\nu_k,\nu_k)
%		=&\nabla_{\nu_{\p \Sigma_{k-1}}}\log {\rho_{k-1}}+H_{\p \Sigma_{k-1}}-H_{\p \Sigma_k}
	\end{align*}
\end{proof}

\begin{proof}[Proof of Theorem \ref{NoFBSlicing}]
	Using stability inequality on the last slice $\Sigma_m$ we have,
	\begin{align*}
		0\leq& \int_{\Sigma_m} -\rho_{m-1}f^2(f^{-1}\Delta f+(\Ric_{\Sigma_{m-1}}(\nu_m)+|\sff_{\Sigma_m}|^2))\\
		&+\int_{\Sigma_m}\rho_{m-1} f^2(\nabla^2\log\rho_{m-1}(\nu_m,\nu_m)-\langle  \nabla_{\Sigma_m}\log \rho_{m-1},\nabla_{\Sigma_m}\log f\rangle-\nabla_{\nu_{m}}h)\\
		&+\int_{\p \Sigma_m} f^2( \nabla_{\nu_{\p \Sigma_m}}{\log f}-\sff_{\p \Sigma_{m-1}}(\nu_m,\nu_m)).
	\end{align*}
	Since $\Sigma_m$ is compact by assumption, we can now plug in $f=\rho_{m-1}^{-1}$ to the above equation and using the same computation for interior terms as in  the proof of Theorem \ref{BM-Cm}, and using the boundary conditions we have,	\begin{align*}
		0\leq &\int_{\Sigma_m} -\rho^{-1}_{m-1}(\Lambda+\cR+\cE+\cG)+\int_{\p \Sigma_m} \rho^{-2}_{m-1}(-\nabla_{\nu_{\p \Sigma_{m-1}}}\log \rho_{m-1}-\sff_{\p \Sigma_{m-1}}(\nu_m,\nu_m))\\
		\leq & \int_{\Sigma_m} -\rho^{-1}_{m-1} (\cC_m(\nu_1,...,\nu_m)+\sum_{k=1}^m \cV_k)
		-\int_{\p \Sigma_m} \rho^{-2}_{m-1} \left(\sum_{l=1}^m\sff_{\p X}(\nu_l,\nu_l)\right)
	\end{align*}
	Under our assumption of $C(m,n)\geq 0$ we have $\cV_k\geq 0$ for $1\leq k\leq m$ as proved in \cite{brendle2024generalization}. We have $\cC_m(\nu_1,...,\nu_m)\geq 0$ by assumption and the boundary terms are positive by the $m$-convexity assumption. This leads to a contradiction.
	
\end{proof}

We now proof that manifolds of the type $(Y,\p Y)\times \TT^m$ admit stable weighted slicing of order $m$.

\begin{proof}[Proof of Theorem \ref{RelHomology}]
	The proof is the same as of Theorem 1.5 in \cite{brendle2024generalization} using relative homology instead.
	
	Let both $N$ and $Y$ be closed oriented smooth manifold with boundary.
	Consider $F:(N,\p N)\rightarrow (X,\p X)= (Y,\p Y)\times T^m$ and $F(\p N) \subset \p X=\p Y \times \TT^m$, a map of non-zero degree. 
	For any top degree form $\omega$ on $X$ that vanishes along $\p X$,
	\begin{equation}
		\int_{N} F^*\omega=\text{deg}(F)\int_{X}\omega.
	\end{equation}
	Now take $f_j$ to be the projection of $F$ onto the $j$-th $\SS^1$-factor in $\TT^m$, and $\omega_j=f^*_j\theta$ with $\theta$ a one form on $\SS^1$ normalized so that it integrates to $1$. Take $f_0: N \rightarrow Y$ the projection map and let $\omega'$  be a top degree form on $Y$ that vanishes along $\p Y$ and normalized so that it integrates to $1$ and denote $\Omega:=f_0^*\omega'$ . We have,
	\begin{equation*}
		\int_N \omega_1\wedge...\wedge\omega_m\wedge \Omega=d.
	\end{equation*}
	
	Now if we are given $\Sigma_{k-1} (k\geq 1)$  with $\int_{\Sigma_{k-1}}\omega_k\wedge...\wedge\omega_m\wedge\Omega=d$, we want to find smooth minimizers $\Sigma_{k}$ to the $\rho_{k-1}$ weighted area functional among the class of $(n-k)$-integer rectifiable currents $\Sigma$ with finite mass  normalized so that $\int_{\Sigma_k} \omega_{k+1} \wedge...\wedge\omega_m\wedge\Omega=d$ and $\p \Sigma_k \subset \p \Sigma_{k-1}\subset ... \p N$. We also have for $\iota_k:\p \Sigma_{k-1}\rightarrow \Sigma_{k-1}$,
	\begin{align*}
		\iota_k^*(\omega_{k+1}\wedge...\wedge\omega_m\wedge\Omega)=0
	\end{align*} 
	
	Now take $p_k$ to be a regular value of $f_k\rvert_{\Sigma_{k-1}}$, and denote $\psi_k$ to be a function of $\SS^1\setminus\{p_k\}$ such that $d\psi_k=\theta$. By the same proof as in \cite{brendle2024generalization} we have $\tilde{\Sigma}_k:=f_k^{-1}(p_k)\cap \Sigma_{k-1}$ is non-empty and,
%	\begin{align*}
%		\omega_k\wedge...\wedge \omega_m\wedge\Omega=d(\psi_k\omega_{k+1}\wedge...\wedge\omega_m\wedge\Omega).
%	\end{align*}
%	Again we have 
by Stokes' theorem,
	\begin{align*}
		d=&\int_{\Sigma_{k-1}\setminus \tilde{\Sigma}_k}\omega_k\wedge...\wedge\omega_m\wedge\Omega\\=&\int_{\Sigma_{k-1}\setminus\tilde{\Sigma}_k} d(\psi_k\omega_{k+1}\wedge...\wedge\omega_m\wedge\Omega)\\
		=&\int_{\p \Sigma_{k-1}} \psi_k\omega_{k+1}\wedge...\wedge\omega_m\wedge\Omega\pm \int_{\tilde{\Sigma}_k}\omega_{k+1}\wedge...\wedge\omega_m\wedge\Omega\\
		=&\pm \int_{\tilde{\Sigma}_k}\omega_{k+1}\wedge...\wedge\omega_m\wedge\Omega
	\end{align*}
	By choosing the orientation of $\tilde{\Sigma}_k$,  we constructed a $\tilde{\Sigma}_k$ with $\p \tilde{\Sigma}_k \subset \p \Sigma_{k-1} $ and $ \int_{\tilde{\Sigma}_k}\omega_{k+1}\wedge...\wedge\omega_m\wedge\Omega=d$. 
	Minimize among this class we obtained free boundary stable solutions to the weighted area functional. Regularity for free boundary area minimizers is studied in \cite{gruter1987optimal}, and no singularity occurs for minimal hypersurfaces in ambient manifold of dimension 7 or less.
	
	In this way we create stable minimal slicing of order $m$ in a manifolds with non-zero degree mapping to $(Y,\p Y)\times \TT^m$ for a closed orientable manifold $Y$, hence cannot have non-negative $m$-intermediate curvature and $m$-convex boundary by Theorem \ref{NoFBSlicing}.
\end{proof}

\begin{proof}[Proof of Theorem \ref{InRadiusTh}]
	We may assume $m\geq 2$. The case $m=1$ was proved in \cite{li2014sharp}. The idea for the case $m=1$ is the following (here $\Sigma_{m-1}=\Sigma_0=X$).
	For any point $p\in \Sigma_{m-1}$, denote $r(p)=\inf_{q\in \p \Sigma_{m-1}} d_{\Sigma_{m-1}}(p,q)<\infty$, then since $\p M \hookrightarrow M$ is proper and $M$ is complete, for any $p\in \Sigma_{m-1}$, we can find $q=q(p)\in \p \Sigma_{m-1}$ such that $d_{\Sigma_{m-1}}(p,q)=r(p)$ is realized by a free boundary minimizing geodesic $l=l(p,q)$.
	One studies the second variation for $l$ which implies the desired bound.

Now let $m\geq 2$, if $\p \Sigma_{m-1}$ is compact, an argument using $\mu$-bubble can be made with the following adaptions from the proof of Theorem \ref{BM-Cm}.

For Caccioppoli sets $\Omega$ with $\Sigma=\p^*\Omega$ in $\Sigma_{m-1}$, we write $$\cA_{m}(\Sigma)=\int_{\Sigma}\rho_{m-1}-\int_{\Sigma_{m-1}}\rho_{m-1}h(\chi_{\Omega}-\chi_{\Omega_0}).$$

	We solve the following minimizing $\mu$-bubble problem for $\cA_{m}$ on the last slice $\Sigma_{m-1}$, take any point $p\in \Sigma_{m-1}$ and any $0<r<r(p):=\inf_{q\in \p \Sigma_{m-1} }d_{\Sigma_{m-1}}(p,q)$, we want to prove $r \leq C_0$.
	Take $h$ to be a smooth mollifier of $\tilde{h}(x):= H_0\frac{r}{r+\epsilon-d_{\Sigma_{m-1}}(x,\p \Sigma_{m-1})}= H_0 \phi(x)$ defined over the set 
	$$\tilde{T}:=\{x\in \Sigma_{m-1},  d_{\Sigma_{m-1}}(x,\p \Sigma_{m-1})<r+\epsilon<r(p)\},$$
	for some small $\epsilon$ and $\Omega_0=\{x\in \Sigma_{m-1}, h(x)>\epsilon\}$ with $\epsilon$ a regular value of $h$.	 We note that by the proof of Lemma 4.1 in \cite{wu2025rigidity}, the following inequality $$h\rvert_{\p \Sigma_{m-1}} -\nabla_{\nu_{\p \Sigma_{m-1}}}\log\rho_{m-1} =H_0\frac{r}{r+\epsilon}-\sum_{i=1}^{m-1}\sff_{\p X}(\nu_i,\nu_i)<H_{\p X}-\sum_{i=1}^{m-1}\sff_{\p X}(\nu_i,\nu_i)=H_{\p \Sigma_{m-1}},$$
	guarantees a non-empty $\Sigma_m=\p^*\Omega_m$ minimizer of $\cA_m$.
	
	Again a minimizing $\mu$-bubble $\Sigma_m \rightarrow \Sigma_{m-1}$ is smooth when the dimension of $\Sigma_{m-1}$ is no more than $7$. Similar to the proof of Theorem \ref{BM-Cm} we obtain the following re-arranged second variation inequality,
	\begin{align*}
		0\leq &-\int_{\Sigma_m} \rho^{-1}_{m-1}(\Lambda+\cR+\cE+\cG +\nabla_{\nu_m}h+H_{\Sigma_m}h)\\
		\leq &-\int_{\Sigma_m}\rho^{-1}_{m-1}(\nabla_{\nu_m}h-\frac{1}{m-1}hH_{\Sigma_m}+C(m,n)H^2_{\Sigma_m}+\frac{m}{2(m-1)}h^2)\\
		0\leq & \int_{\Sigma_{m}}-\rho^{-1}_{m-1}(C_0 h^2+\nabla_{\nu_m}h)\\
		\leq &\int_{\Sigma_m}-\rho_{m-1}^{-1}(C_0H_0\phi^2(x)-\phi^2(x)\cdot \frac{1}{r}),
		\end{align*}
and we obtained $C_0 H_0 r\leq 1$.
\end{proof}

\section{Spectral Non-negative Ricci}\label{BMSpec}

\begin{defi}[Spectral Ricci Curvature]\label{AlmSpecRic}
	For a constant $\gamma\geq 0$, we say that a complete Riemannian manifold $M$ has almost non-negative Ricci curvature in the spectral sense if,
	 there is $\alpha' \in (0,1)$ and $u\in C^{2,\alpha'}(M)$ such that $u>0$ and 
	 \begin{equation*}
	 	-\gamma \Delta u+\Ric \cdot u =\lambda_1 u-\alpha u|\nabla \log u|^2,
	 \end{equation*}
	 for some $\gamma>0, \alpha\geq 0, \lambda_1 \geq 0$ and $\lambda_1=0$ when $M$ is non-compact.

	Here we denote $\Ric(x)$ as $\Ric_M(x):=\inf_{v \in T_xM, g(v,v)=1} \Ric(v,v) \in \text{Lip}_{\text{loc}}(M)$.
\end{defi}

\begin{defi}[Spectral Mean Convexity]\label{SpecMC}
	For a constant $\gamma\geq 0$, we say that a complete Riemannian manifold $M$ has non-negative  mean curvature in the spectral sense if for some $H_0 \geq 0$,
	\begin{align*}
		\gamma \nabla_{\nu_{\p M}} \log u+H_{\p M}\geq H_0.
	\end{align*}
\end{defi}

\begin{lemm}
	If $X$ is a complete Riemannian manifold with non-negative $m$-intermediate curvature $\cC_m \geq 0$ and non-negative mean curvature $H_{\p X} \geq H_0\geq 0$, and $X=\Sigma_0\supset \Sigma_1 ... \supset \Sigma_{m-1}$ is a stable weighted free boundary slicing of order $m-1$. Assume $m^2-2-n(m-2)\geq 0$, then $\Sigma_{m-1}$ has almost non-negative Ricci curvature and non-negative mean curvature in the spectral sense for $\gamma=1$ and $\alpha=(\frac{1}{2}-\frac{1}{2(m-1)})$.
\end{lemm}

\begin{proof}
	By Lemma \ref{BHJ1} and Lemma \ref{FB-BC} we have that over $\Sigma_{m-1}$,
	\begin{align*}
		-\Delta_{\Sigma_{m-1}}\log \rho_{m-1}&= \Lambda+ \cE'+\cR' +\cG,\\
		\nabla_{\nu_{\p \Sigma_{m-1}}}\log u&=\sum_{i=1}^{m-1}\sff_{\p X}(\nu_i,\nu_i)=H_{\p X}-H_{\p \Sigma_{m-1}}\geq H_0- H_{\p \Sigma_{m-1}}
	\end{align*}
	 here $\cE'=\sum_{k=1}^{m-1}|\sff_{\Sigma_k}|^2-H^2_{\Sigma_k}$, $\cR'=\sum_{i=1}^{m-1}\Ric_{\Sigma_{k-1}}(\nu_k)$, and $\Lambda, \cG$ as defined in the proof of Theorem \ref{BM-Cm}.
	 
	 Let $\Ric(x):=\Ric_{\Sigma_{m-1}}(x)$ be obtained at $\Ric(\nu_m)$ for some $\nu_m \in \Sigma_{m-1}$ and by definition $\nu_{m} \perp \{\nu_1,...,\nu_{m-1}\}$, then denoting $\alpha_{m-1}=\frac{m-2}{2(m-1)}=\frac{1}{2}-\frac{1}{2(m-1)},$
	 \begin{align*}
	 	-\Delta_{\Sigma_{m-1}}\log \rho_{m-1}+\Ric=& \Lambda+\cG +\cE'+\cR\\
	 	\geq & \cC_m(\nu_1,...,\nu_{m-1},\nu_m)+\sum_{k=1}^{m-1} \cV_k+(1-\alpha_{m-1})|\nabla_{\Sigma_{m-1}}\log \rho_{m-1}|^2,
	 \end{align*}
	 here $\cV_k \geq 0$ for $1\leq k\leq m-1$ when $m^2-2-n(m-2)\geq 0$.
	 Now using,
	 \begin{align*}
	 	-\Delta_{\Sigma_{m-1}} \log \rho_{m-1}=-\rho^{-1}_{m-1}\Delta \rho_{m-1}+|\nabla_{\Sigma_{m-1}} \log \rho_{m-1}|^2
	 \end{align*}
	 we obtain,
	 \begin{align*}
	 	-\Delta_{\Sigma_{m-1}}\rho_{m-1}+\Ric \cdot \rho_{m-1} \geq -\alpha_{m-1}\rho_{m-1}|\nabla_{\Sigma_{m-1}}\log \rho_{m-1}|^2.
	 \end{align*}
\end{proof}

\begin{lemm}[Second Variation Formula]\label{SecVarSpec}
	Assume $\Sigma^k$ is a complete smooth stable critical points of the following functional on $M^{k+1}$, for some $\gamma\geq 0$ and $h$ a smooth function,
	\begin{align*}
		\cA_f(\Omega)=\int_{\Sigma}f^{\gamma}-\int_{M}(\chi_{\Omega}-\chi_{\Omega_0})hf^{\gamma}
	\end{align*}
	The first variation implies $H_{\Sigma}=h-\gamma\nabla_{\nu_{\Sigma}}\log f$.
	
 The second variation  gives the stability inequality, for any $\phi \in C^1_c(\Sigma)$,
	 \begin{align*}
	 	0\leq & \int_{\Sigma}f^{\gamma}(-\phi\Delta\phi-\phi^2(|\sff_{\Sigma}|^2+\Ric_M(\nu_{\Sigma})))+\gamma f^{\gamma-1}\phi^2 (\Delta_M f-\Delta_{\Sigma}f-H_{\Sigma}\nabla_{\nu_{\Sigma}}f)\\
		&+\int_{\Sigma} -f^{\gamma}\nabla_{\nu_{\Sigma}} h\phi^2-\gamma f^{\gamma}\phi^2 |\nabla_{\nu_{\Sigma}}\log f|^2-\gamma\phi f^{\gamma-1}\nabla_{\Sigma}f\cdot \nabla_{\Sigma}\phi
	 \end{align*}
\end{lemm}

\begin{theo}\label{TheoSpec}
	Assume $M^{k+1} (k+1\leq 7)$ has almost positive Ricci curvature in the spectral sense, that is for some $k_0>0$ and $\gamma>0, \alpha\geq 0$, we have a positive smooth function $f$ over $M$ with,
	\begin{align*}
	 	-\gamma\Delta_{{M}}f+\Ric \cdot f \geq k_0f -\alpha f|\nabla_{M}\log f|^2.
	 \end{align*}
	 Then 
	 \begin{enumerate}
	 	\item when $k=1$, $\alpha<1,(\gamma-2)^2<4(1-\alpha)$, then $\text{diam}(M^2)\leq \frac{\pi}{\sqrt{c_1k_0}}$ for $c_1=1-\frac{\gamma^2}{4(\gamma-\alpha)}>0$. 
	 	\item when $k=2$, if either $\alpha=0,0<\gamma\leq 2$ or $0<\alpha<\frac{1}{2}, 1-\sqrt{1-2\alpha}<\gamma <1+\sqrt{1-2\alpha}$ then $\text{diam}(M^3)\leq \frac{\pi}{\sqrt{\frac{k_0}{2}}}$.
	 	\item when $k>2, \alpha<\frac{1}{k}$, $\frac{2}{k}(1-\sqrt{1-k\alpha})<\gamma<\frac{2}{k}(1+\sqrt{1-k\alpha})$, then $\text{diam}(M^{k+1})\leq \frac{\pi}{\sqrt{C(\alpha,k,\gamma)k_0}}$ for  $C(\alpha,k,\gamma)=\frac{-k\gamma^2+4\gamma-4\alpha}{4k(\frac{\gamma^2}{k}-\gamma^2+\gamma-\alpha)}>0$.
	 \end{enumerate}
\end{theo}

\begin{proof}[Proof of Theorem \ref{TheoSpec}] Given a $\mu$-bubble solution $\Sigma$ to $\cA_f(\cdot)$ as in Lemma \ref{SecVarSpec} for a suitable function $h$, we have from the second variation in Lemma \ref{SecVarSpec},
	\begin{align*}
		0\leq &\int_{\Sigma}f^{\gamma}|\nabla_{\Sigma}\phi|^2-\phi^2(|\sff_{\Sigma}|^2+\Ric_M(\nu_{\Sigma}))f^{\gamma}+\gamma f^{\gamma-1}\phi^2 (\Delta_M f-\Delta_{\Sigma}f-H_{\Sigma}\nabla_{\nu_{\Sigma}}f)\\
		&+\int_{\Sigma} -f^{\gamma}\nabla_{\nu_{\Sigma}} h\phi^2-\gamma f^{\gamma}\phi^2 |\nabla_{\nu_{\Sigma}}\log f|^2
	\end{align*}
	We plug in $\phi=\varphi f^{\frac{-\gamma}{2}}$ and use $H_{\Sigma}=h-\gamma\nabla_{\nu_{\Sigma}}\log f$, $|\sff_{\Sigma}|^2 \geq \frac{1}{k}H_{\Sigma}^2$ to obtain,
	\begin{align*}
		0\leq &\int_{\Sigma}|\nabla_{\Sigma}\varphi|^2-k_0\varphi^2-\frac{1}{k}\varphi^2h^2-(\nabla_{\nu_{\Sigma}}h) \varphi^2+|\nabla_{\Sigma} \log f|^2\varphi^2 (\frac{\gamma^2}{4}-\gamma+\alpha)\\
		&+\int_{\Sigma} \varphi^2 (\nabla_{\nu_{\Sigma}}\log f)^2(\gamma^2-\gamma-\frac{\gamma^2}{k}+\alpha)+\varphi^2 h\nabla_{\nu_{\Sigma}}\log f (\frac{2}{k}-1)\gamma+\gamma\varphi\nabla \varphi\cdot \nabla\log f.
	\end{align*}
	Since
	
	When $k=2$, and assume $\frac{\gamma^2}{4}-\gamma+\alpha\leq 0$ and $\gamma^2-\gamma-\frac{\gamma^2}{k}+\alpha\leq 0$ we obtain
	\begin{align*}
		0\leq \int_{\Sigma}-k_0-\frac{1}{2}h^2+|\nabla_{\nu_{\Sigma}}h|.
	\end{align*}
	The two conditions together gives $0<\gamma\leq 2$ when $\alpha=0$, which coincides with condition (1.5) in \cite{xu2025dimension}. When $\alpha>0$, the two conditions simplifies to $\alpha\leq \frac{1}{2}$ and $1-\sqrt{1-2\alpha}< \gamma <1+\sqrt{1-2\alpha}$. 
	Arguing as in the proof of Theorem \ref{BM-Cm}, we obtain $\text{diam}(M^3)\leq \frac{\pi}{\sqrt{\frac{k_0}{2}}}$.

	Now assume $k\neq 2$ and $\frac{\gamma^2}{4}-\gamma+\alpha\leq 0$ and $\gamma^2-\gamma-\frac{\gamma^2}{k}+\alpha< 0$, we have by Young's inequality,
	\begin{align*}
		0\leq &\int_{\Sigma} C'(\alpha,k,\gamma)|\nabla_{\Sigma}\varphi|^2+\varphi^2\left(\frac{(\frac{2}{k}-1)^2\gamma^2}{4(\frac{\gamma^2}{k}+\gamma-\gamma^2-\alpha)}-\frac{1}{k} \right)h^2+\varphi^2 |\nabla_{\nu_{\Sigma}}h|-k_0\varphi^2\\
		=&\int_{\Sigma}C'(\alpha,k,\gamma)|\nabla_{\Sigma}\varphi|^2 -\varphi^2(k_0+ C(\alpha,k,\gamma)h^2-|\nabla_{\nu_{\Sigma}}h|).
	\end{align*}
	Arguing as in the proof of Theorem \ref{BM-Cm}, we obtain $\text{diam}(M)\leq \frac{\pi}{\sqrt{C(\alpha,k,\gamma)k_0}}$ when $C(\alpha,k,\gamma)>0$.

	We compute for $\gamma>0$,
	\begin{align*}C(\alpha,k,\gamma)=\frac{1}{k}-\frac{(\frac{2}{k}-1)^2\gamma^2}{4(\frac{\gamma^2}{k}+\gamma-\gamma^2-\alpha)}=\frac{4\gamma-4\alpha-k\gamma^2}{4k(\frac{\gamma^2}{k}+\gamma-\gamma^2-\alpha)}>0
	\end{align*}
	has a solution exactly when $\alpha<\frac{1}{k}$ and,
	\begin{align}\label{Spec1}
	\frac{2}{k}(1-\sqrt{1-k\alpha})<\gamma<\frac{2}{k}(1+\sqrt{1-k\alpha}) \iff \frac{k}{4}\gamma^2-\gamma+\alpha<0.
	\end{align}
	When $k=1$, the two conditions $\frac{\gamma^2}{4}-\gamma+\alpha\leq 0$ and $\gamma^2-\gamma-\frac{\gamma^2}{k}+\alpha<0$ gives,
	\begin{align*}
		\alpha <\gamma-\frac{\gamma^2}{4},
	\end{align*}
	which is equivalent to (\ref{Spec1}) for $k=1$.
	
	When $k>2$, the inequality (\ref{Spec1}) implies   $\gamma^2-\gamma-\frac{\gamma^2}{k}+\alpha<0$, which implies $\frac{\gamma^2}{4}-\gamma+\alpha\leq 0$.

%	
%	This gives when $\alpha<\frac{k}{4(k-1)}$,
%	\begin{align}
%		\frac{k}{2(k-1)}(1- \sqrt{1-4\alpha(1-\frac{1}{k})})<\gamma< \frac{k}{2(k-1)}(1+ \sqrt{1-4\alpha(1-\frac{1}{k})}). \label{Spec3}
%	\end{align}
	When $\alpha=0$, the three conditions give $0<\gamma<\frac{4}{k}$ and coinside with the condition (1.5) in \cite{xu2025dimension}.
\end{proof}

\begin{theo}\label{TheoSpecB}
	Assume $(M^{k+1}, \p M), k+1\leq 7$ has almost positive Ricci curvature and uniformly positive mean curvature in the spectral sense, for some smooth positive function $f$ and constant $k_0>0$, $\gamma>0, \alpha\geq 0$ and $H_0>0$ as in Definition \ref{AlmSpecRic}, Definition \ref{SpecMC}.	If $\p M$ is compact, then we have the following bound on the in-radius of $M$, $\text{Rad}_{in}(M):=\sup_{x\in M}d_{M}(x,\p M)$,
	 \begin{enumerate}
	 	\item when $k=1$, $\alpha<1,(\gamma-2)^2<4(1-\alpha)$, then $\text{Rad}_{in}(M^2)\leq \frac{1}{{c_1H_0}}$ for $c_1=1-\frac{\gamma^2}{4(\gamma-\alpha)}>0$. 
	 	\item when $k=2$, if either $\alpha=0,0<\gamma\leq 2$ or $0<\alpha<\frac{1}{2}, 1-\sqrt{1-2\alpha}<\gamma <1+\sqrt{1-2\alpha}$ then $\text{Rad}_{in}(M^3)\leq \frac{2}{H_0}$.
	 	\item when $k>2, \alpha<\frac{1}{k}$, $\frac{2}{k}(1-\sqrt{1-k\alpha})<\gamma<\frac{2}{k}(1+\sqrt{1-k\alpha})$, then $\text{Rad}_{in}(M^{k+1})\leq \frac{1}{{C(\alpha,k,\gamma)H_0}}$ for  $C(\alpha,k,\gamma)=\frac{-k\gamma^2+4\gamma-4\alpha}{4k(\frac{\gamma^2}{k}-\gamma^2+\gamma-\alpha)}>0$.
	 \end{enumerate}
\end{theo}

\begin{proof}[Proof of Theorem \ref{TheoSpecB}] 
	Similar to the proof of Theorem \ref{InRadiusTh}, we would like to find a $\mu$-bubble solution to $\cA_f(\cdot)$ as in Lemma \ref{SecVarSpec} for a suitable $h$ such that $h\rvert_{\p M}<H_0$.
	
	By a first variation argument similar to Lemma 4.1 in \cite{wu2025rigidity}, the spectral mean-convexity condition,
	\begin{align*}
		h\rvert_{\p M}-\gamma\nabla_{\p M}\log f<H_0-\nabla_{\nu_{\p M}}\log f\leq H_{\p M}
	\end{align*}
	guarantees a non-empty solution to $\cA_f(\cdot)$.
	
	The rest of the argument follows exactly as in the proof of Theorem \ref{TheoSpec}.
\end{proof}

\bibliographystyle{alpha}
\bibliography{REF.bib}

\begin{thebibliography}{MWY25}

\bibitem[BHJ24]{brendle2024generalization}
Simon Brendle, Sven Hirsch, and Florian Johne.
\newblock A generalization of geroch's conjecture.
\newblock {\em Communications on Pure and Applied Mathematics}, 77(1):441--456, 2024.

\bibitem[Che24]{chen2024generalization}
Shuli Chen.
\newblock A generalization of the geroch conjecture with arbitrary ends.
\newblock {\em Mathematische Annalen}, 389(1):489--513, 2024.

\bibitem[CKL25]{chu2025rigidity}
Jianchun Chu, Kwok-Kun Kwong, and Man-Chun Lee.
\newblock Rigidity on non-negative intermediate curvature.
\newblock {\em Mathematical Research Letters}, 31(6):1693--1714, 2025.

\bibitem[CL24]{chodosh2024generalized}
Otis Chodosh and Chao Li.
\newblock Generalized soap bubbles and the topology of manifolds with positive scalar curvature.
\newblock {\em Annals of Mathematics}, 199(2):707--740, 2024.

\bibitem[GL83]{gromov1983positive}
Mikhael Gromov and H~Blaine Lawson.
\newblock Positive scalar curvature and the dirac operator on complete riemannian manifolds.
\newblock {\em Publications Math{\'e}matiques de l'IH{\'E}S}, 58:83--196, 1983.

\bibitem[Gro19]{gromov2019mean}
Misha Gromov.
\newblock Mean curvature in the light of scalar curvature.
\newblock In {\em Annales de l'Institut Fourier}, volume~69, pages 3169--3194, 2019.

\bibitem[Gr{\"u}87]{gruter1987optimal}
Michael Gr{\"u}ter.
\newblock Optimal regularity for codimension one minimal surfaces with a free boundary.
\newblock {\em manuscripta mathematica}, 58(3):295--343, 1987.

\bibitem[Li14]{li2014sharp}
Martin Man-chun Li.
\newblock A sharp comparison theorem for compact manifolds with mean convex boundary.
\newblock {\em The Journal of Geometric Analysis}, 24(3):1490--1496, 2014.

\bibitem[Mag12]{maggi2012sets}
Francesco Maggi.
\newblock {\em Sets of finite perimeter and geometric variational problems: an introduction to Geometric Measure Theory}, volume 135.
\newblock Cambridge University Press, 2012.

\bibitem[MWY25]{mazurowski2025topology}
Liam Mazurowski, Tongrui Wang, and Xuan Yao.
\newblock On the topology of manifolds with positive intermediate curvature.
\newblock {\em arXiv preprint arXiv:2503.13815}, 2025.

\bibitem[SY79]{yau1979structure}
R~Schoen and ST~Yau.
\newblock On the structure of manifolds with positive scalar curvature.
\newblock {\em Manuscripta mathematica}, 28:159--184, 1979.

\bibitem[SY96]{shen1996stable}
Ying Shen and Rugang Ye.
\newblock On stable minimal surfaces in manifolds of positive bi-{Ricci} curvatures.
\newblock {\em Duke Math. J.}, 85(1):109--116, 1996.

\bibitem[Wu23]{wu2023free}
Yujie Wu.
\newblock Free boundary stable minimal hypersurfaces in positively curved 4-manifolds.
\newblock {\em arXiv preprint arXiv:2308.08103}, 2023.

\bibitem[Wu25]{wu2025rigidity}
Yujie Wu.
\newblock Rigidity of complete free boundary minimal hypersurfaces in convex nnsc manifolds.
\newblock {\em arXiv preprint arXiv:2504.20585}, 2025.

\bibitem[Xu25]{xu2025dimension}
Kai Xu.
\newblock Dimension constraints in some problems involving intermediate curvature.
\newblock {\em Transactions of the American Mathematical Society}, 378(03):2091--2112, 2025.

\end{thebibliography}

\end{document}